\documentclass[12pt]{amsart}
\usepackage{amsfonts}
\usepackage{graphicx,amssymb,amsfonts,amsmath,amscd}
\usepackage[all,ps,cmtip]{xy}
\usepackage{amscd}

\usepackage{mathrsfs}
\let\mathcal\mathscr

\usepackage{hyperref}

\def\dra{\dashrightarrow}

\def\Alb{\mathop{\rm Alb}\nolimits}

\def\codim{\mathop{\rm codim}\nolimits}
\def\dim{\mathop{\rm dim}\nolimits}

\def\Ker{\mathop{\rm Ker}\nolimits}

\def\llra{\hbox to 12mm{\tofill}}

\def\Pic{\mathop{\rm Pic}\nolimits}

\def\Bs{\mathop{\rm Bs}\nolimits}

\def\Q{{\bf Q}}

\def\a{{\alpha}}

\def\t{{\tau}}

\def\cF{{\mathcal F}}

\def\cH{{\mathcal H}}
\def\cI{{\mathcal I}}
\def\cJ{{\mathcal J}}

\def\cL{{\mathcal L}}

\def\cO{{\mathcal O}}

\def\cQ{{\mathcal Q}}

\def\1Y{{Y'}}
\def\aY{{\widetilde{Y}}}
\def\aX{{\widetilde{X}}}

\def\WW{{\widehat{W}}}

\def\ax{{\overline{X}}}
\def\ay{{\overline{Y}}}
\def\av{{\overline{V}}}

\def\af{{\overline{f}}}

\newtheorem{theo}{Theorem}[section]

\newtheorem{lemm}[theo]{Lemma}

\newtheorem{rema}[theo]{Remark}

\newtheorem{exam}[theo]{Example}

\begin{document}
\title{On varieties of maximal Albanese dimension}
\date{\today}
\author{Zhi Jiang}
\maketitle
A smooth projective complex variety $X$ has {\em maximal Albanese dimension} if its Albanese
  map $ X\to \Alb(X)$ is generically finite onto its image. These varieties have recently attracted a lot of attention and have been shown to have very special geometric properties
(\cite{cha},  \cite{chb}, \cite{asia}, \cite{fuj}, \cite{HAC1},   \cite{par},  \cite{pp}).

Assume for example that $f: X\to Y$ is a surjective  morphism between smooth projective varieties  of the same dimension. For each   positive integer $m$, denote by $P_m(X):=h^0(X,\omega_X^{\otimes m})$ the $m$-th plurigenus of $X$. We have $P_m(X)\ge P_m(Y)$, but it is in general difficult
  to conclude anything on $f$ if there is  equality.
However, when $Y$ (hence also $X$) is of general type and has maximal Albanese
dimension, Hacon and Pardini proved  in \cite{HAC1}, Theorem 3.2, that if  $P_m(X)=P_m(Y)$ for some $m\geq 2$,  then $f$ is
birational. We give in   \S\ref{exa}  examples that show that this conclusion does not hold in general.

More generally, if $X\to I(X)$ and $Y\to I(Y)$ are the respective Iitaka
fibrations of $X$ and  $Y$, we may assume, taking appropriate birational models,  that $f$ induces a morphism $I(f): I(X)\to I(Y)$. When $Y$ has maximal Albanese
dimension, but is not necessarily of general type, Hacon and Pardini proved  that  if  $P_m(X)=P_m(Y)$ for some $m\geq 2$, then $I(f)$ has
connected fibers (since $I(Y)$ is birational to $Y$ when $Y$ is of general type, this implies the result quoted above).

But in their proof, Hacon and Pardini actually do not use the
assumption that $Y$ has maximal Albanese dimension; all they
need is that $P_m(X)=P_m(Y)>0$ and $I(Y)$ has maximal Albanese
dimension (see section 1). However, under their assumption, we prove here a much stronger
conclusion.

\medskip
\noindent{\bf Theorem 1} {\em Let $f: X\to Y$ be a surjective  morphism between smooth complex projective varieties  of the same dimension. If $Y$ has maximal Albanese dimension and  $P_m(X)=P_m(Y)$ for some $m\geq 2$, the induced map $I(f): I(X)\dra I(Y)$ between the respective Iitaka models of $X$ and $Y$
  is birational.

  Moreover, $f$ is birationally equivalent to a quotient by a finite abelian group.}
  \medskip

For more details on the last statement, we refer to
Theorem \ref{thm7}.
\medskip

In another direction, it was shown by Chen and Hacon  (\cite{cha}, Theorem 4) that if $X$ is a smooth projective
 variety of  maximal Albanese dimension, the image of the 6-canonical map $\phi_{6K_X}$ has dimension the Kodaira dimension $\kappa(X)$. If $X$ is moreover of general type, $\phi_{6K_X}$ is  birational onto its image (\cite{asia}, Corollary 4.3). We prove a common generalization of these results (Theorem \ref{iit}):

\medskip
\noindent{\bf Theorem 2} {\em  If $X$ is a smooth complex projective variety with maximal Albanese dimension, $\phi_{5K_X}$ is a model of the Iitaka fibration of $X$.}
\medskip

The proof follows the ideas of \cite{pp2} and is based on a result from \cite{jiang}. We also prove that
$\phi_{3K_X}$ is already  a model of the Iitaka fibration of $X$ under a stronger assumption on $X$  (Theorem \ref{iit2}).
\medskip

The article is organized as follows. The first section is devoted to the proof of the birationality of $I(f)$. In the second section, we give a complete structure theorem  for $f$ (Theorem \ref{thm7}) which shows that the situation is quite restricted. In the third  section, we present three examples showing that the conclusion of the above theorem   can fail when the varieties do not have maximal Albanese dimension, and in the last section, we prove our results on pluricanonical maps of varieties of   maximal Albanese dimension.

We work over the field of complex numbers.
\section{Proof that $I(f)$ is birational}

 We begin with a general lemma (we refer to \cite{Laz}, \S11, for the definition and properties of the  asymptotic multiplier ideal sheaf $\cJ(||D||)$ associated with a divisor $D$ on a smooth projective variety).

\begin{lemm}\label{6} Let $f: X\to Y$ be a surjective morphism
between smooth projective varieties of the same dimension with $\kappa(Y)\geq 0$. For
any $m\geq 2$,
$$f_*(\cO_X(mK_X)\otimes\cJ(||(m-1)K_X||))\supset \cO_Y(mK_Y)\otimes\cJ(||(m-1)K_Y||).$$
\end{lemm}
\begin{proof}Take $N>0$  and let $\t_{Y}: Y'\to Y$ be a log-resolution
such that $$\t_Y^*|N(m-1)K_Y|=|L_1|+E_1,$$ where $|L_1|$ is
base-point-free and $E_1$ is the fixed divisor. Then we take  a log-resolution $\t_X: X'\to X$ such that we have a
commutative diagram:
$$\CD
  X' @>f'>> Y' \\
  @V \t_X VV @V \t_Y VV  \\
  X @>f>> Y
\endCD$$
and $\t_X^*|N( m-1)K_X|=|L_2|+E_2$ where $|L_2|$ is base-point-free
and $E_2$ is the fixed divisor. Let $D\in |(m-1)K_{X/Y}|$. Then
${f'}^* E_1+N\t_X^*D\succeq E_2$. Hence
\begin{eqnarray*}&&\cO_{X'}(K_{X'/X}+m\t_X^*K_X-\Big\lfloor\frac{1}{N}E_2 \Big\rfloor)\\&\supset&
\cO_{X'}(K_{X'/X}+m\t_X^*K_X-\t_X^*D-\Big\lfloor\frac{1}{N}{f'}^*E_1 \Big\rfloor)\\
&=&\cO_{X'}(K_{X'/X}+\t_X^*K_X+(m-1)\t_X^*f^*K_Y-\Big\lfloor\frac{1}{N}{f'}^*E_1 \Big\rfloor)\\
&=&\cO_{X'} \big(K_{X'/Y'}-\Big\lfloor\frac{1}{N}{f'}^*E_1 \Big\rfloor+{f'}^* \Big\lfloor\frac{1}{N}E_1 \Big\rfloor+{f'}^*(K_{Y'/Y}+m\t_Y^*K_Y-\Big\lfloor\frac{1}{N}E_1 \Big\rfloor) \big).
\end{eqnarray*}
We may assume that $N$ is sufficiently large and divisible. Then
$$\t_{X*} \big(\cO_{X'}(K_{X'/X}+m\t_X^*K_X-\Big\lfloor\frac{1}{N}E_2 \Big\rfloor) \big)=\cO_X(mK_X)\otimes\cJ(||(m-1)K_X||).$$
By step 2 in the proof of
\cite{HAC1}, Theorem 3.2, we know that
$K_{X'/Y'}-\big\lfloor\frac{1}{N}{f'}^*E_1 \big\rfloor+{f'}^* \big\lfloor\frac{1}{N}E_1 \big\rfloor$
is an effective divisor, hence
\begin{eqnarray*}
&& \t_{Y*}f'_* \Big(\cO_{X'} \big(K_{X'/Y'}-\Big\lfloor\frac{1}{N}{f'}^*E_1 \Big\rfloor +{f'}^* \Big\lfloor\frac{1}{N}E_1 \Big\rfloor\\
&&\hskip 4cm {}+
{f'}^*(K_{Y'/Y}+m\t_Y^*K_Y-\Big\lfloor\frac{1}{N}E_1 \Big\rfloor) \big) \Big)\\
& \supseteq&
\t_{Y*} \big(\cO_{Y'}(K_{Y'/Y}+m\t_Y^*K_Y-\Big\lfloor\frac{1}{N}E_1 \Big\rfloor) \big)\\
 &=& \cO_Y(mK_Y)\otimes\cJ(||(m-1)K_Y||).\end{eqnarray*}
This proves the lemma.
\end{proof}

 We now prove the first part of Theorem 1, stated in the introduction. We start from
  a surjective  morphism  $f: X\to Y$ between smooth projective varieties  of the same dimension.

  Changing the notation from the introduction, we let $  V$ and $ W$ be the respective Iitaka
models of $X$ and  $Y$, and
we may assume, taking appropriate birational models, that we have a commutative
diagram of {\em morphisms}
\begin{equation}\label{23}
\xymatrix{
  X\ar[r]^{f}\ar[d]^{h_X}&Y\ar[d]^{h_Y}\ar[r]^{a_Y}& A\ar[d]^{\pi}\\
  V\ar[r]^{g}& W\ar[r]^{a_W} &A/K
}
\end{equation}
where  $h_X$ and $h_Y$ are the respective Iitaka fibrations  of $X$ and $Y$,
 $a_Y$ and $a_W$ are the respective Albanese morphisms of $Y$ and $W$, and $K$ is an abelian subvariety of $A:=\Alb(Y)$
(see \cite{HAC1}, \S2.1).
 We set
\begin{eqnarray*}
 \cH_X:=h_{X*}(\cO_X(mK_X)\otimes\cJ(||(m-1)K_X||))&\quad&\cF_X:=a_{W*}g_*\cH_X
\\
\cH_Y:=h_{Y*}(\cO_Y(mK_Y)\otimes\cJ(||(m-1)K_Y||))& \quad&\cF_Y:=a_{W*} \cH_Y.
\end{eqnarray*}
When $m\geq 2$,
we have
$
\cF_Y\subset\cF_X$ by   Lemma \ref{6} and we denote by  $\cQ$ the quotient sheaf  $\cF_X/
\cF_Y$ on $A/K$.

{\em Assume now}   $P_m(X)=P_m(Y)=M>0$. By Theorem 11.1.8 and Proposition 11.2.10 in
\cite{Laz}, we have
\begin{eqnarray*}
 P_m(Y)&=&h^0(Y,
\cO_Y(mK_Y)\otimes\cJ(||mK_Y||))\\
&=&h^0(Y,
\cO_Y(mK_Y)\otimes\cJ(||(m-1)K_Y||))\\
&=&h^0(W,
\cH_Y)\\
&=&h^0(A/K, \cF_Y).
\end{eqnarray*}
Similarly,
$$
 P_m(X)=h^0(V, \cH_X)=h^0(A/K, \cF_X).
$$
Thus $\cH_Y\subset h_{Y*}(\cO_Y(mK_Y))$ is a nonzero torsion-free
sheaf. Since $h_Y$ is a model of the Iitaka fibration of $Y$ whose
general fibers are birationally isomorphic  to abelian varieties (\cite{HAC1}, Proposition
2.1), the latter sheaf has rank $1$. So the rank of
$\cH_Y$ is also $1$. We have the same situation for $h_X$, hence the
rank of
$\cH_X$ is again $1$.
On the other hand, we claim the following.

\medskip{\bf Claim:} $\cQ=0$, hence
$\cF_Y=\cF_X$.

In order to prove the Claim, we want to apply Proposition 2.3 in
\cite{HAC1}. Namely, it is enough to prove   $h^j(A/K,
\cF_Y\otimes P)=h^j(A/K, \cF_X\otimes P)$  for all $j\geq 0$ and all
$P\in\Pic^0(A/K).$
 We will first prove that when $j\geq 1$. By Lemma 2.1 in \cite{jiang},
 we have
\begin{equation}\label{26}
H^i(W, \cH_Y\otimes a_W^*P)=H^i(W, g_*\cH_X\otimes a_W^*P)=0,
\end{equation}
for all $P\in \Pic^0(A/K)$ and all $i\geq 1$.
 We now prove
\begin{equation}\label{27}
R^ja_{W*}\cH_Y=R^ja_{W*}(g_*\cH_X)=0,
\end{equation}
for all $j\geq 1$, as follows.

First we take a very ample line
bundle $H$ on $A/K$ such that, for all $k\geq 1$
and $j\geq 0$,
 \begin{equation}\label{28}H^k(A/K, R^ja_{W*}\cF_Y\otimes H)=H^k(A/K, R^ja_{W*}(g_*\cF_X)\otimes H)=0\end{equation}
  and $R^ja_{W*}\cF_Y\otimes H$ and $R^ja_{W*}(g_*\cF_X)\otimes H$ are globally generated. Again by Lemma 2.1 in \cite{jiang}, $$H^j(W, \cH_Y\otimes a_W^*H)=H^j(W, g_*\cH_X\otimes a_W^*H)=0,$$ for all
$j\geq 1$. Therefore, by Leray's spectral sequence and (\ref{28}), we conclude that
$$H^0(A/K, R^ja_{W*}\cH_Y\otimes H)=H^0(A/K, R^ja_{W*}(g_*\cH_X)\otimes H)=0
$$
 for $j\geq 1$. Since $R^ja_{W*}\cH_Y\otimes H$ and $R^ja_{W*}(g_*\cH_X)\otimes H$ are globally generated, we deduce that $R^ja_{W*}\cH_Y=R^ja_{W*}(g_*\cH_X)=0$, for all $j\geq 1$.

Applying the Leray spectral sequence to (\ref{26}), we get,  by (\ref{27}),
for all $i\geq 1$ and $P\in \Pic^0(A/K)$,
$$H^i(A/K, \cF_Y\otimes P)=H^i(W, \cH_Y\otimes g^*P)=0,
$$
 and
$$H^i(A/K, \cF_X\otimes P)=H^i(W, g_*\cH_X\otimes g^*P)=0.$$

Finally, for all $P\in\Pic^0(A/K)$,
\begin{eqnarray*}h^0(A/K,
\cF_Y\otimes P)&=&\chi(A/K, \cF_Y\otimes P)=\chi(A/K,
\cF_Y)\\&=&h^0(A/K, \cF_Y)=M,\end{eqnarray*}
and similarly,
$$h^0(A/K, \cF_X\otimes P)=h^0(A/K, \cF_X)=M.$$

We have finished the proof of the Claim.

\medskip
{\em Assume moreover that $W$ has maximal Albanese dimension,}   so that
  $a_W$ is generically
finite onto its image $Z$, the rank of $\cH_Y$ is 1, and the rank of
$\cF_X=\cF_Y=a_{W*}\cH_Y$ on $Z$ is $\deg(a_W)$.
Consider the Stein factorization
$$g:V\xrightarrow{p}U\xrightarrow{q}W,$$
where $p$ is an algebraic fiber space and $q$ is surjective and
finite. Because $h^0(U, p_*\cH_X)=h^0(V, \cH_X)=M>0$, the nonzero torsion-free sheaf
$p_*\cH_X$ has rank $\geq 1$. We can
write
 $$\cF_X=a_{W*}g_*\cH_X=a_{W*}q_*(p_*\cH_X) $$
  and conclude that the rank of $\cF_X$ on $Z$ is $\geq \deg(q)\cdot\deg(a_W)$.
 This implies $\deg(q)=1$ hence $g$ has connected fibers. Essentially, this is  Hacon and Pardini's proof of \cite{HAC1}, Theorem 3.2.

\medskip
{\em Assume finally that $Y$ has maximal Albanese dimension.}  We just saw that $g\circ h_X$ is an
algebraic fiber space  and we denote by $X_w$ a general fiber. The main ingredient is the following lemma.

\begin{lemm}\label{ma}In the above situation, the sheaf $$g_*\cH_X=(g\circ h_X)_*(\cO_X(mK_X)\otimes\cJ(||(m-1)K_X||))$$ has
rank $P_m(X_w)>0$.
\end{lemm}

This lemma will be proved later. We first use it to finish the proof
of the first part of Theorem 1.

Assume that $g$
is not birational.  Since it is an algebraic
fiber space, we have $\dim (W)<\dim (V)$. Hence by the easy addition formula
(\cite{Mo}, Corollary 1.7), we have $\dim (V)=\kappa(X)\leq
\kappa(X_w)+\dim (W)$, hence $\kappa(X_w)\geq 1$. Since $X$ is of
maximal Albanese dimension, $X_w$ is also of maximal Albanese
dimension, hence $P_{m}(X_w)\geq 2$ by Chen and Hacon's
characterization of abelian varieties (\cite{CH1}, Theorem 3.2).
Then, by Lemma \ref{ma}, the rank of $\cH_X$ on $Z$ is
$\deg(a_W)\cdot P_m(X_w)$ ($\geq 2\deg(a_W)$), which is a
contradiction. This concludes the proof.

\medskip
In order to prove Lemma \ref{ma}, we begin with an easy lemma.

\begin{lemm}\label{5}Let $X$ be a smooth projective variety, let $D_1$ be a divisor
on $X$ with nonnegative Iitaka dimension, and let $D_2$  be an effective
divisor on $X$. We have an inclusion
$$\cJ(||D_1+D_2||)\supset
\cJ(||D_1||)\otimes\cO_X(-D_2).$$
\end{lemm}

\begin{proof}Take $N>0$ such that
$|ND_1|\neq\emptyset$. Choose a log-resolution $$\mu:
X'\to X$$ for $ND_1$, $ND_2$, and $N(D_1+D_2)$. Write
\begin{eqnarray*}
\mu^*(|ND_1|)&=&|W_1|+E_1\\
\mu^*(|ND_2|)&=&|W_2|+E_2\\
\mu^*(|N(D_1+D_2)|)&=&|W_3|+E_3,
\end{eqnarray*}
where $E_1$, $E_2$, and $E_3$ are the fixed divisors and $|W_1|$,
$|W_2|$, and $|W_3|$ are free linear series. We have
\begin{gather*}N\mu^*D_2\succeq E_2\qquad \textmd{and} \qquad E_1+E_2\succeq E_3,
\end{gather*}
hence
\begin{eqnarray*}
\mu_*(K_{X'/X}-\Big\lfloor\frac{1}{N}E_3 \Big\rfloor)&\supset&
\mu_*(K_{X'/X}-\Big\lfloor\frac{1}{N}(E_1+E_2) \Big\rfloor)\\
&\supset&\mu_*(K_{X'/X}-\Big\lfloor\frac{1}{N}(E_1+N\mu^*D_2) \Big\rfloor)\\
&=&\mu_*(K_{X'/X}-\Big\lfloor\frac{1}{N}E_1 \Big\rfloor)\otimes\cO_X(-D_2).
\end{eqnarray*}
By the definition of asymptotic multiplier ideal sheaves, this
proves Lemma \ref{5}.
\end{proof}

\begin{proof}[Proof of Lemma \ref{ma}]
We will reduce Lemma \ref{ma} to Proposition 3.6 in \cite{jiang}. Since $Y$
is of maximal Albanese dimension and $h_Y$ is a model of the Iitaka
fibration of $Y$, by a theorem of Kawamata (see also Theorem 3.2 in \cite{jiang}), there exists an \'{e}tale cover $\pi_Y:
\widetilde{Y}\to Y$ induced by an \'{e}tale cover of $A$ and
a commutative diagram:
$$
\xymatrix{
\widetilde{Y}\ar[r]^{\pi_Y}\ar@{-->}[d]_{h_{\aY}}& Y\ar[d]^{h_Y}\\
\WW\ar[r]^{b_{\WW}}&W,}
$$
where $\WW$ is a smooth projective variety of general type, the
rational map $h_{\aY}$ is a model of the Iitaka fibration of $\aY$,
and $b_{\WW}$ is generically finite and surjective.

Let $\aX$ be a connected component of $X\times_Y\widetilde{Y}$, denote by
$\pi_{\aX}$ the induced morphism $\aX\to X$, and
denote by $f_{\aX}$ the induced morphism $\aX\to\aY$. Denote
by $k$ and $k_{\aX}$ respectively the morphism $g\circ h_X=h_Y\circ
f$ and the map $h_{\aY}\circ f_{\aX}$. After birational
modifications of $\aX$, we may suppose that $k_{\aX}$ is a morphism such that $k_{\aX}(E)$ is a proper
subvariety of $\WW$, where $E$ is the $\pi_{\aX}$-exceptional divisor.
All in all, we have the commutative diagram:
$$
\xymatrix@R=30pt@C=30pt@M=+5pt{\aX\ar[r]^{\pi_{\aX}}\ar[d]^{f_{\aX}}\ar@/_2pc/[dd]_{k_{\aX}}&X\ar[d]^f\ar@/^2pc/[dd]^{k}\\
\widetilde{Y}\ar[r]^{\pi_Y}\ar@{-->}[d]^{h_{\aY}}& Y\ar[d]^{h_Y}\\
\WW\ar[r]^{b_{\WW}}&W .}
$$
We then take the Stein factorization:
$$k_{\aX}:\aX\xrightarrow{k_1} W_1\xrightarrow{b_{W_1}}\WW.$$
The important point is that $W_1$ is still of general type. Again by taking
birational modifications of $\aX$ and $W_1$, we may assume that
$k_1: \aX\to W_1$ is an algebraic fiber space between smooth
projective varieties. We can apply Proposition 3.6 in \cite{jiang} to the following diagram:
$$
\xymatrix{
\widetilde{X}\ar[rr]^{\pi_{\aX}}\ar[d]^{k_1}&&X\ar[d]^k\\
W_1\ar[rr]_{b_{\WW}\circ b_{W_1}}&&W.}
$$
It follows that the sheaf
$$k_*(\cO_X(mK_X)\otimes\cJ(||(m-1)K_{X/W}+k^*K_W||))\otimes\cO_W(-(m-2)K_W)$$ has rank
$P_m(X_w)$. By Lemma 3.4 in \cite{jiang}, the line bundle
$(m-1)K_{X/W}+k^*K_W$ has nonnegative Iitaka dimension. By Lemma
\ref{5},
$$\cJ(||(m-1)K_X||)\supset \cJ(||(m-1)K_{X/W}+k^*K_W||)\otimes
\cO_X(-(m-2)k^*K_{W}).$$ Therefore,
\begin{eqnarray*}
&&k_*\big(\cO_X(mK_X) \big)\\
&\supset& k_* \big(\cO_X(mK_X)\otimes\cJ(||(m-1)K_X||) \big)
\\
&\supset& k_* \big(\cO_X(mK_X)\otimes\cJ(||(m-1)K_{X/W}+k^*K_W||) \big)\otimes\cO_W \big(-(m-2)K_W \big).\end{eqnarray*}
Since the rank of the first and the third sheaves are both $P_m(X_w)$, so is the rank of the
second.
\end{proof}

\section{A complete description of $f: X\to Y$}

By using Kawamata's Theorem 13 in \cite{KA} (see also
Theorem 3.2 in \cite{jiang}), we obtain the following complete description of $f$.

 \begin{theo}\label{thm7}
Let $f: X\to Y$ be a surjective morphism of smooth
  projective varieties of the same dimension, with $Y$ of maximal Albanese
dimension.

If $P_m(X)=P_m(Y)$ for some
$m\geq 2$, there exist  \begin{itemize}
\item a normal projective variety $V_X$ of general type,
\item an abelian variety $A_X$,
\item
a finite abelian group $G$ which acts faithfully on $V_X$ and on
$A_X$ by translations,
\item a subgroup $G_2$ of $G$,
\end{itemize}
such that
\begin{itemize}
\item   $X$ is
birational to $(A_X\times V_X)/G$, where $G$  acts diagonally on $A_X\times V_X$,
\item   $Y$ is
birational to  $(A_Y\times V_Y)/G_1$, where
\subitem $\star$  $V_Y=V_X/G_2$ and $A_Y=A_X/G_2$,
\subitem $\star$ $G_1:=G/G_2$ acts diagonally on
$A_Y\times V_Y$,
\item  $f$ is
birational to the quotient morphism $(A_X\times V_X)/G\to
(A_Y\times V_Y)/G_1$.
\end{itemize}
\end{theo}

\begin{proof}
In the diagram (\ref{23}), we already know that $g: V\to W$ is
birational so we may assume that $V=W$ and $g$ is the identity. We
then consider the diagram:
\begin{equation}\label{29}
\xymatrix{
X\ar[r]^{f}\ar[d]^{h_X}&Y\ar[r]^{a_Y}\ar[d]^{h_Y}&A\ar[d]\\
V\ar@{=}[r]& V\ar[r]^{a_V} & A/K. }
\end{equation}
Taking the Stein factorizations for $f$ and $a_Y$, we may assume that $X$ and $Y$ are normal and $f$ and $a_Y$ are finite. Similarly we take the Stein factorization for $Y\xrightarrow{a_Y} A\to A/K$ and may assume that $V$ is normal and $a_V$ is finite.

By Poincar\'{e} reducibility, there exists an isogeny $B\to
A/K$ such that $A\times_{A/K}B\simeq K\times B$. We denote by $H$
the kernel of this isogeny. Apply the \'{e}tale base change
$B\to A/K$ to diagram (\ref{29}) and get
\begin{equation}\label{30}
\xymatrix{
\ax\ar[r]^{\af}\ar[d]^{h_{\ax}}&\ay\ar[r]^-{a_{\ay}}\ar[d]^{h_{\ay}}& K\times B\ar[d]\\
\av\ar@{=}[r]& \av\ar[r]^{a_{\av}} & B, }
\end{equation}
where
\begin{itemize}
\item $\av=V\times_{A/K}B$ and $\ay=Y\times_V\av$ (which are connected because $a_{Y}$ and $a_V$ are the Albanese maps),
\item
  $\ax=X\times_Y\ay$ (which is also connected because $\ax=X\times_Y\ay=X\times_Y(Y\times_V\av)=X\times_V\av$),
  \item
$h_X: X\to V$ is an algebraic fiber space.
\end{itemize}
Let $A_X$ and $A_Y$ be the respective general fibers of $h_{\ax}$ and $h_{\ay}$. We have the following induced diagram from (\ref{30}):
$$
\xymatrix{ A_X\ar[r]_{\beta}\ar@/^1pc/[rr]^{\alpha_X}& A_Y\ar[r]_{\alpha_Y}& K
}
$$
 By Proposition 2.1 in \cite{HAC1}, $A_X$ and $A_Y$
are birational to abelian varieties. Hence the morphisms $\alpha_X$
and $\alpha_Y$ are birationally equivalent to \'{e}tale covers.
Since $a_{\ay}$ and $\a_{\ay}\circ \af$ are finite, $\alpha_X$ and
$\alpha_Y$ are also finite. Thus $\alpha_X$ and $\alpha_Y$ are
isogenies of abelian varieties by Zariski's Main Theorem. We denote
by $\widetilde{G}$, $\widetilde{G}_1$, and $\widetilde{G}_2$ the
abelian groups $\Ker(A_X\to K)$, $\Ker(A_Y\to K)$, and $\Ker(A_X\to A_Y)$ respectively. Then
$\widetilde{G}_1=\widetilde{G}/\widetilde{G}_2$ and
$A_Y=A_X/\widetilde{G}_2$. Let $k\in K$ be a general point, let
$V_Y$ be the normal variety $a_{\ay}^{-1}(k\times B)$, and let $V_X$ be
the normal variety $\af^{-1}a_{\ay}^{-1}(k\times B)$.

We know that $A_X$ and $A_Y$ respectively act on $\ax$ and $\ay$ in
such a way that $\af$ is equivariant for the $A_X$-action on $\ax$
and the $A_Y$-action on $\ay$.
Furthermore, the actions induce a faithful $\overline{G}$-action on
$V_X$ and a faithful $\widetilde{G}_1$-action on $V_Y$, and we have
an $A_X$-equivariant isomorphism $\ax\simeq(A_X\times
V_X)/\widetilde{G}$ and an $A_Y$-equivariant isomorphism $\ay\simeq
(A_Y\times V_Y)/\widetilde{G}_1$, where $\widetilde{G}$ acts on
$A_X\times V_X$ diagonally and $\widetilde{G}_1$ acts on $A_Y\times
V_Y$ diagonally.

The  induced morphism
$$
\xymatrix{
V_X\ar[r]^{\af\vert_{V_X}}\ar[d]^{h_{\ax}}& V_Y\ar[d]^{h_{\ay}}\\
V\ar@{=}[r]&V }
$$
  is equivariant for the $\widetilde{G}$-action on $V_X$ and the $\widetilde{G}_1$-action on $V_Y$. Thus $V_Y=V_X/\widetilde{G}_2$ and $\af\vert_{V_X}$ is the quotient morphism.

Thus we obtain  $A_Y=A_X/\widetilde{G_2}$ and $V_Y=V_X/\widetilde{G_2}$, and $\af: \ax\to\ay$ is the quotient morphism $(A_X\times V_X)/\widetilde{G}\to (A_Y\times V_Y)/\widetilde{G}_1$, so $$\af: \ax=(A_X\times V_X)/\widetilde{G}\to \ay=(A_Y\times V_Y)/\widetilde{G}_1$$ is also the quotient morphism.

Let $G=\Ker(A_X\times B\to A)$ and $G_1=\Ker(A_Y\times B\to A)$. We have
exact sequences of groups
$$
1\to \widetilde{G}\to G\to H\to 1
\quad \rm{and}\; \quad 1\to \widetilde{G}_1\to
G_1\to H\to 1.
$$
Then $X=(A_X\times V_X)/G$ and $Y=(A_Y\times V_Y)/G_1$, and $f$ is
the quotient map. This proves Theorem \ref{thm7} with $G_2= \widetilde{G}_2\subset G$.

\end{proof}

 \section{Examples}\label{exa}

In the next two examples, we see that the conclusion of our theorem does not hold in general, even for surfaces of general type.

\begin{exam}\label{e1}\upshape Let $C_1$ and $C_2$ be smooth projective curves of genus $2$ with
respective hyperelliptic involutions
$i_1$ and $i_2$. Define $Y$ to be the minimal resolution of singularities of $(C_1\times C_2)/(i_1, i_2)$. Let $X$ be the blow-up of $C_1\times C_2$ at the $36$ fixed points of $(i_1, i_2)$. There is a 2-to-1 morphism $f: X\to Y$. We have $K_Y^2=\frac{1}{2}K_{C\times C}^2=4$ and $c_2(Y)=\frac{1}{2}(c_2(C_1\times C_2)-36)+72=56$. Since $Y$ is a minimal surface, we have $P_2(Y)=K_Y^2+\frac{1}{12}(K_Y^2+c_2(Y))=9$. We also have $P_2(X)=P_2(C_1\times C_2)=3\times 3=9$. Hence we have a nonbirational morphism $f: X\to Y$ between smooth projective surfaces of general type with $q(Y)=0$ (so that $Y$ does not have maximal Albanese dimension!) and $P_2(X)=P_2(Y)=9$.
\end{exam}

\begin{rema}\upshape It turns out that the situation in the case of surfaces of general type can be completely worked out. More precisely, one can show that
if $f: S\rightarrow T$ is a nonbirational morphism between smooth projective surfaces
of general type such that $P_m(S)=P_m(T)$ for some $m\geq 2$, then $m=2$ and one of the following occurs:
\begin{itemize}
\item[1)] either $S$ is birational to the product of two smooth projective  curves of genus 2, and $f$ is birationally equivalent to the quotient by the diagonal hyperelliptic involution (see Example \ref{e1} above);
\item[2)] or $S$ is birational to the   theta divisor  of the Jacobian of a smooth projective  curve  of genus 3 and $f$ is birationally equivalent to the bicanonical map of $S$;
\item[3)] $S$ is birational to a double cover of a principally polarized abelian surface  branched along a divisor in $ |2\Theta|$ having at most double points and $f$ is birationally equivalent to the bicanonical map of $S$.
\end{itemize}
\end{rema}

In higher dimensions, we have many more examples.

\begin{exam} \label{e2}\upshape(Compare with \cite{K}, Proposition 8.6.1)
We denote by $\mathbf{P}(a_0^{s_0}, \ldots, a_k^{s_k})$  the weighted projective space with $s_j$ coordinates of
weight $a_i$ (see \cite{D}). For any integer $k\geq 3$, denote by $P_X$ the weighted projective
space $\mathbf{P}(1, (2k)^{4k+5}, (2k+1)^{4k-3})$ with coordinates $x_i$ and
by $P_Y$ the weighted projective space $\mathbf{P}(2, (2k)^{4k+5}, (2k+1)^{4k-3})$ with coordinates $y_i$.
As in the proof of Proposition 8.6.1 in \cite{K}, one can check that $P_X$ and $P_Y$
both  have canonical singularities. There is a natural degree-2
morphism $\varepsilon: P_X\to P_Y$ defined by $y_0=x_0^2$ and $y_i=x_i$ for $i\geq 1$.

Let $Y'$ be a general hypersurface of weighted degree $d=16k^2+8k$ in $P_Y$ and let $X'$
be the pull-back by $\varepsilon$ of $Y'$. Since $2k(2k+1)|d$ and $Y'$ is general,
$X'$ is also general and both $X'$ and $Y'$ have canonical singularities. Take
resolutions $X\to X'$ and $Y\to Y'$ such that $\varepsilon$ induces  a degree-2 morphism
$f: X\to Y$. The canonical sheaves are $\omega_{X'}=\cO_{X'}(2)$ and $\omega_{Y'}=\cO_{Y'}(1)$.
Since both $X'$ and $Y'$ have canonical singularities, we have, for any integer $m\geq 0$,
$$P_m(X)=h^0(X', \cO_{X'}(2m))  \quad \textrm{and}\quad P_m(Y)=h^0(Y', \cO_{Y'}(m)).$$
It follows from Theorem in 1.4.1 in \cite{D} that for $m$ even and $<2k$, we have
$P_m(X)=P_m(Y)=1$. By Theorem 4.2.2 and Corollary 2.3.6 in \cite{D}, $q(X)=q(Y)=0$.
\end{exam}

Under the assumptions of our theorem, one might
expect that $f: X\to Y$   be birational to an \'{e}tale
morphism. However the example below (see also Example 1 in \cite{CH4})
shows that
this is not the case in general.

\begin{exam}\upshape Let $G=\mathbb{Z}_{rs}$ and let $G_2=s\mathbb{Z}_{rs}$ be the subgroup of $G$ generated by $s$, with $s\geq 2$ and  $r\geq
2$. Let $G_1=G/G_2\simeq \mathbb{Z}_s$. Consider an elliptic curve
$E$, let $B_1$ and $B_2$ be two points on $E$, and let $L$ be
a line bundle of degree $1$ such that $B=(rs-a)B_1+aB_2\in |tmL|$
with $1\leq a\leq m-2$ and $(a, rs)=1$. Taking the normalization of the
$(rs)$-th root of $B$, we get a smooth curve $C$ and a Galois cover
$\pi: C\to E$ with Galois group $G$. By construction, $\pi$
ramifies at two points, $B_1$ and $B_2$. Following \cite{Be} \S VI.12, we have $h^0(C,
\omega_C^2)^G=2$.

Let $L^{(i)}$ be $L^i(-\lfloor\frac{iB}{rs}\rfloor)$ and denote
$(L^{(i)})^{-1}$ by $L^{-(i)}$. Then, by Proposition 9.8 in \cite{K1},
$$\pi_*\cO_C=\bigoplus_{i=0}^{rs-1}L^{-(i)}.$$ Let $C_1$ be the curve
$$\underline{Spec}(\bigoplus_{i=0}^{s-1}L^{-(ri)}),$$ where
$\bigoplus_{i=0}^{s-1}L^{-(ri)}$ has the subalgebra structure of
$\pi_*\cO_C$. Consider the Stein factorization
$$\pi:C\xrightarrow{g} C_1\xrightarrow{\pi_1} E.$$ Then $C_1=C/G_2$ and
$\pi_1$ is a Galois cover with Galois group $G_1$ which also
ramifies only at $B_1$ and $B_2$. Hence we again have
$$h^0(C_1,
\omega_{C_1}^2)^{G_1}=2.$$

Finally we take an abelian variety $K$ such that $G$ acts freely on
$K$ by translations and set $K_1=K/G_2$. Let
$$\widetilde{X}=C\times
K\quad{\rm ,}\quad \widetilde{Y}=C_1\times K$$ and
$$X=\widetilde{X}/G=(C\times
K)/G\quad{\rm ,}\quad Y=\widetilde{Y}/G=(C_1\times K_1)/G_1,$$ where $G$ and $G_1$
act diagonally. Hence $\widetilde{X}$ and $\widetilde{Y}$ are
\'{e}tale covers of $X$ and $Y$ respectively. There is a natural
finite dominant morphism $f: X\to Y$ of degree $r$. Since
its lift $\widetilde{f}: \widetilde{X}\to \widetilde{Y}$ is
not \'{e}tale, $f$ is not \'{e}tale.

Since $$H^0(X, \omega_X^2)\simeq H^0(\widetilde{X}, \omega_{\widetilde{X}}^2)^G \simeq  H^0(C, \omega_C^2)^G$$ and
$$H^0(Y, \omega_Y^2)\simeq H^0(\widetilde{Y}, \omega_{\widetilde{Y}}^2)^{G_1} \simeq  H^0(C_1, \omega_{C_1}^2)^{G_1},$$
 we have
$$P_2(X)=P_2(Y)=2.$$
\end{exam}

\section{Pluricanonical maps of varieties of maximal Albanese dimension}\label{last}

Let $X$ be a smooth projective variety of maximal Albanese dimension. As mentioned in the introduction,
Chen and Hacon proved in \cite{cha} and \cite{asia} that
$\phi_{6K_X}(X)$ has dimension $\kappa(X)$; if $X$ is moreover of general type,
$\phi_{6K_X}$ is  birational onto its image. They also showed that if $\chi(\omega_X)>0$,
the map $\phi_{3K_X}$ is already birational onto its image.
Pareschi and Popa provided in \cite{pp2}, \S6,  a conceptual approach to
these theorems based on their regularity and vanishing theorems.

We prove a unifying statement for  varieties with maximal Albanese dimension which are not necessarily of general type.
The proof
is parallel to that of Pareschi and Popa. 

In this section, we will always assume $f: X\rightarrow Y$ is a birational model of the Iitaka fibration of $X$.

\begin{theo}\label{iit}
If $X$ is a smooth projective variety with maximal Albanese dimension, the linear system $|\cO_X(5K_X)\otimes f^*P|$ induces
the Iitaka fibration. In particular,
$\phi_{5K_X}$ is a model of the Iitaka fibration of $X$.
\end{theo}

\begin{proof} We may as in (\ref{23}) assume that we have a diagram
\begin{eqnarray}
\xymatrix{ X\ar[d]_f\ar[r]^-{a_X}&\Alb(X)\ar[d]^{f_*}\\
Y\ar[r]^-{a_Y}&\Alb(Y)}\end{eqnarray}
where $f$ is the Iitaka fibration of $X$ and
 $a_X$ and $a_Y$ are the respective Albanese morphisms of $X$ and $Y$.

Since $f$ is a model of the Iitaka fibration of $X$, $f_*(\omega_X^2\otimes\cJ(||K_X||))$
is a torsion-free rank $1$ sheaf on $Y$. We now use the following lemma (\cite{jiang},  Lemma 2.1):

\begin{lemm}\label{mul}Suppose that $f: X\rightarrow Y$ is a surjective
morphism between smooth projective varieties, $L$ is a
$\Q$-divisor on $X$, and the Iitaka model of $(X, L)$
dominates $Y$. Assume that $D$ is a nef $\Q$-divisor on $Y$
such that $L+f^*D$ is a divisor on $X$. Then we have$$H^i(Y,
R^jf_*(\cO_X(K_X+L+f^*D)\otimes \cJ(||L||)\otimes Q))=0,$$ for all
$i\geq 1$, $j\geq 0$, and all $Q\in \Pic^0(X)$.
\end{lemm}
By Lemma \ref{mul}, we have
$$H^i(Y, f_*(\cO_X(2K_X)\otimes\cJ(||K_X||)\otimes Q))=0$$
for all $i\geq 1$ and $Q\in\Pic^0(X)$. As in Lemma 2.6 in \cite{jiang},
 $R^ja_{Y*}(f_*(\cO_X(2K_X)\otimes\cJ(||K_X||)\otimes Q))=0$ for all $j\geq 1$. Hence
\begin{eqnarray*}
&&H^i(\Alb(Y), a_{Y*}f_*(\cO_X(2K_X)\otimes\cJ(||K_X||)\otimes Q))\\
&=&H^i(Y, f_*(\cO_X(2K_X)\otimes\cJ(||K_X||)\otimes Q))\\
&=&0,
\end{eqnarray*}
for all $i\geq 1$ and $Q\in\Pic^0(X)$. Thus for any $Q\in V_0(\omega_X)\subset \Pic^0(X)$,
$a_{Y*}f_*(\cO_X(2K_X)\otimes\cJ(||K_X||)\otimes Q)$ is a nonzero IT-sheaf of index $0$ and in particular,
it is $M$-regular. By \cite{pp2}, Corollary 5.3,  $a_{Y*}f_*(\cO_X(2K_X)\otimes\cJ(||K_X||)\otimes Q)$ is
continuously globally generated. Since $a_Y$ is generically finite,
the exceptional locus $Z_1$ of $a_Y$ is a proper closed subset of $Y$.
Then $f_*(\cO_X(2K_X)\otimes\cJ(||K_X||)\otimes Q)$ is continuously globally generated away from $Z_1$.
By definition, this means that for any open subset $V\subset\Pic^0(Y)$, the evaluation map
\begin{multline*}\bigoplus_{P\in V}H^0(Y, f_*(\cO_X(2K_X)\otimes\cJ(||K_X||)\otimes Q)\otimes P)\otimes P^{-1}
\\\rightarrow f_*(\cO_X(2K_X)\otimes\cJ(||K_X||)\otimes Q)\end{multline*}
 is surjective away from $Z_1$.

Now we claim that there exists an open dense subset $U\subset Y-Z_1$ such that the sheaf $$a_{Y*}(\cI_y\otimes f_*(\cO_X(3K_X)\otimes \cJ(||2K_X||)))$$ is $M$-regular for any $y\in U$.

We first assume the claim and finish the proof of the theorem.

We conclude by the claim that $\cI_y\otimes f_*(\cO_X(3K_X)\otimes \cJ(||2K_X||))$ is continuously 
globally generated away from $Z_1$. Denote respectively by $\cL$ and $\cL_1$ the rank-$1$ torsion-free sheaves
$f_*(\cO_X(2K_X)\otimes\cJ(||K_X||))$ and $f_*(\cO_X(3K_X)\otimes\cJ(||2K_X||))$. Let $U_1$ be a dense open subset 
of $Y-Z_1$ such that $\cL$ and $\cL_1$ are locally free on $U_1$.
Then by \cite{pp3}, Proposition 2.12, $\cL_1\otimes \cL$ is very ample over $U\cap U_1$.
We have $\cL\otimes\cL_1\hookrightarrow f_*(\cO_X(5K_X))$, thus $f_*(\cO_X(5K_X))$
is very ample on a dense open subset of $Y$. This concludes the proof of the theorem.

For the claim, let $$U\subset U_1 \bigcap
\big(Y-\bigcup_{T_i}\bigcap_{Q\in T_i}\Bs(|f_*(\omega_X\otimes Q)|)\big)$$ be any dense open subset of $Y$,
where $T_i$ runs through all the components of $V_0(\omega_X)$ and $\Bs(|f_*(\omega_X\otimes Q)|)$ denotes the locus
where the evaluation map $$H^0(Y, f_*(\omega_X\otimes Q))\otimes\cO_Y\rightarrow f_*(\omega_X\otimes Q)$$ is not surjective.
For each component $T_i$ of $V_0(\omega_X)$, we may write $T_i=P_i+f^*S_i$, where $S_i$ is a subtorus of $\Pic^0(Y)$ and
$P_i\in\Pic^0(X)$ (see \cite{GL2}).

Again, by Lemma \ref{mul}, we have $$H^i(Y, \cL_1\otimes Q)=0$$
 for all $i\geq 1$ and any $Q\in \Pic^0(Y)$. For $y\in U$,
consider the exact sequence
 $$0\rightarrow \cI_y\otimes \cL_1\rightarrow \cL_1\rightarrow \mathbb{C}_y\rightarrow 0.$$ We push forward this short sequence
 to $\Alb(Y)$. Since
$y\in U$, we have
$$0\rightarrow a_{Y*}(\cI_y\otimes \cL_1)\rightarrow a_{Y*}\cL_1\rightarrow \mathbb{C}_{a_Y(y)}\rightarrow 0.$$

Hence $H^i(Y, a_{Y*}(\cI_y\otimes \cL_1)\otimes Q)=0$ for any $i\geq 2$ and $Q\in \Pic^0(Y)$.
We now assume that $a_{Y*}(\cI_y\otimes \cL_1)$ is not $M$-regular. Then by definition of $M$-regularity, we have
\begin{eqnarray}\label{nMregular}\codim_{\Pic^0(Y)}V_1(a_{Y*}(\cI_y\otimes \cL_1))\leq 1.\end{eqnarray}

Hence $y$ is a base-point of all sections in $H^0(Y, \cL_1\otimes P_s)$, for all
$s\in V_1(a_{Y*}(\cI_y\otimes \cL_1))$.

On the other hand, by \cite[Lemma 2.2]{jiang}, $$\dim H^0(X, \cO_X(3K_X)\otimes f^*P)$$
is constant for $P\in\Pic^0(Y)$. Then, \begin{eqnarray*}\dim H^0(Y, \cL_1\otimes P)&=&\dim H^0(Y, \cL_1)=
\dim H^0(X, \cO_X(3K_X))\\&=&\dim H^0(X, \cO_X(3K_X)\otimes f^*P).\end{eqnarray*}
Hence the inclusion
$$H^0(X, \cO_X(3K_X)\otimes\cJ(||2K_X||)\otimes f^*P)\hookrightarrow H^0(X, \cO_X(3K_X)\otimes f^*P)$$ is an isomorphism.
Therefore, $y\in \Bs(|f_*\cO_X(3K_X)\otimes P_s|)$, for
all $s\in V_1(a_{Y*}(\cI_y\otimes \cL_1))$.

Since $y\in Y-\bigcup_{T_i}\bigcap_{Q\in T_i}\Bs(|f_*(\omega_X\otimes Q)|)$, let $V_i\subset S_i$ be a dense open subset
such that $y\notin \Bs|f_*(\omega_X\otimes Q)|$, for any $Q\in P_i+f^*V_i$.

We may shrink $U$ so that $f_*(\cO_X(K_X)\otimes P_i)$ and $f_*(\cO_X(2K_X)\otimes P_i^{-1})$
are locally free on $U$ for all $i$. Moreover, we can require that, for each $i$, the multiplication
$$f_*(\cO_X(K_X)\otimes P_i)\otimes f_*(\cO_X(2K_X)\otimes P_i^{-1})\rightarrow f_*\cO_X(3K_X)$$
is an isomorphism on $U$, since both sheaves are of rank $1$.

 We then conclude that $y$ is a base point of all sections of
$$H^0(Y, f_*(\cO_X(2K_X)\otimes P_i^{-1})\otimes Q^{'})$$
where $Q^{'}\in V_1(a_{Y*}(\cI_y\otimes \cL_1))-V_i$.

We may further shrink $U$ so that
$$f_*(\cO_X(2K_X)\otimes \cJ(||K_X||)\otimes P_i^{-1})|_U=f_*(\cO_X(2K_X)\otimes P_i^{-1})|_U$$
is locally free for each $i$. Then $y\in U$ belongs to
$$\Bs|f_*(\cO_X(2K_X)\otimes \cJ(||K_X||)\otimes P_i^{-1})\otimes Q^{'}|$$
for each $Q^{'}\in V_1(a_{Y*}(\cI_y\otimes \cL_1))-V_i$.

By \cite{cha}, Theorem 1, the union of all the $S_i$ generates $\Pic^0(Y)$.
Hence by (\ref{nMregular}), for some $i$, $V_1(a_{Y*}(\cI_y\otimes \cL_1))-V_i$ contains an open subset of $\Pic^0(Y)$ and this contradicts
the fact that $f_*(\cO_X(2K_X)\otimes \cJ(||K_X||)\otimes P_i^{-1})$ is continuous globally generated away from $Z_1$. 
This concludes the proof of the claim.
\end{proof}
Our Theorem \ref{iit} is just an analog of Theorem 6.7 in \cite{pp2}. The main point is just that $a_{Y*}f_*(\cO_X(2K_X)\otimes\cJ(||K_X||))$ is $M$-regular. On the other hand, if $X$ is of general type, of maximal Albanese dimension, and if moreover $a_X(X)$
is not ruled by tori, Pareschi and Popa proved that $a_{X*}\omega_X$ is $M$-regular, which is the main ingredient of the proof of Theorem 6.1 in \cite{pp2}. If $X$ is not of general type,  $a_X(X)$   is always ruled by tori of dimension $n-\kappa(X)$. But
we still have:

\begin{theo}\label{iit2}If $X$ is a smooth projective variety with maximal Albanese dimension $n$, and if
its Albanese image $a_X(X)$ is not ruled by tori of dimension $>n-\kappa(X)$, the map
$\phi_{3K_X}$ is a model of the Iitaka fibration of $X$.
\end{theo}

\begin{proof}We just need to show that under our assumptions, and with the notation of the proof of Theorem \ref{iit}, $a_{Y*}f_*(\omega_X)$ is $M$-regular. The rest is the same as the proof of Theorem \ref{iit}. By Kawamata's theorem \cite[Theorem 13]{KA}, we have the following commutative diagram:

\begin{eqnarray*}
\xymatrix{
\widehat{Y}\times \widetilde{K}\ar[d]^{pr_1}& \widetilde{X}\ar[l]_(.4){\mu}\ar[r]^-{\pi_X}\ar[d]^{\widehat{f}}&X\ar[r]^(.4){a_X}\ar[d]^f& \Alb(X)\ar[d]^{f_*}\\
\widehat{Y}\ar@{=}[r]&\widehat{Y}\ar[r]^{b_Y}&Y\ar[r]^-{a_Y}&\Alb(Y),}
\end{eqnarray*}
where $\pi_X$ is birationally equivalent to a finite \'{e}tale cover of $X$ induced by isogeny of $\Alb(X)$, $\mu$ is a birational morphism, $\widetilde{K}$ is an abelian variety isogenous to $\ker f_*$, $\widehat{Y}$ is a smooth projective variety of general type, and $b_Y$ is generically finite. We   set $g_Y=a_Y\circ b_Y$.

Since $a_X(X)$ is not ruled by tori of dimension $>n-\kappa(X)$,
we conclude that $g_Y(\widehat{Y})=a_Y(Y)$ is not ruled by tori. We make the following:

\smallskip\noindent{\bf Claim:} $g_{Y*}\omega_{\widehat{Y}}$ is $M$-regular.

We first see how the Claim implies Theorem \ref{iit2}. Since $\widetilde{K}$ is an abelian variety,
we have obviously $pr_{1*}\omega_{\widehat{Y}\times\widetilde{K}}=\omega_{\widehat{Y}}$.
Hence $$g_{Y*}pr_{1*}\omega_{\widehat{Y}\times\widetilde{K}}=g_{Y*}pr_{1*}\mu_*\omega_{\widetilde{X}}=
a_{Y*}f_*\pi_{X*}\omega_{\widetilde{X}}$$
is $M$-regular on $\Alb(Y)$. On the other hand, $\omega_X$ is a direct summand of
$\pi_{X*}\omega_{\widetilde{X}}$ since $\pi_X$ is birationally equivalent to an \'{e}tale cover. Therefore, $a_{Y*}f_*\omega_X$ is a direct summand of $g_{Y*}pr_{1*}\omega_{\widehat{Y}\times\widetilde{K}}$ and hence is $M$-regular.

We now prove the Claim.

We first define the following subset of $\Pic^0(Y)$ for any $i\geq 0$:
$$V_i(\widehat{Y}, \Pic^0(Y)):=\{P\in \Pic^0(Y): H^i(\widehat{Y}, \omega_{\widehat{Y}}\otimes g_Y^*P)\neq 0\}.$$
Since the image of $g_Y: \widehat{Y}\rightarrow \Alb(Y)$
is not ruled by tori,
the same argument in the last part of the proof of Theorem 3 in \cite{EL} shows that $\codim_{\Pic^0(Y)}V_i(\widehat{Y}, \Pic^0(Y))>i$ for any $i\geq 1$.
On the other hand, by Grauert-Riemenschneider vanishing, $R^ig_{Y*}\omega_{\widehat{Y}}=0$
for any $i\neq 0$. Thus $$H^i(\widehat{Y}, \omega_{\widehat{Y}}\otimes g_Y^*P)\simeq H^i(\Alb(Y), g_{Y*}\omega_{\widehat{Y}}\otimes P).$$
Hence we have $V_i(g_{Y*}\omega_{\widehat{Y}})=V_i(\widehat{Y}, \Pic^0(Y))$ as subset of $\Pic^0(Y)$. This finishes the proof of the Claim.
\end{proof}

\end{document}